\documentstyle[12pt,leqno,amssymb,amsfonts,amscd,graphics,hyperref,tikz]{amsart}  
\hypersetup{hidelinks=true}

\newtheorem{thm}{Theorem}[section]
\newtheorem{lemma}[thm]{Lemma}
\newtheorem{prop}[thm]{Proposition}
\newtheorem{cor}[thm]{Corollary}
\theoremstyle{definition}
\newtheorem{dfn}[thm]{Definition}
\theoremstyle{remark}

\begin{document}

\newcommand{\y}{Y^\alpha}
\newcommand{\al}{(\alpha_0, \alpha_1, \alpha_2)}
\newcommand{\ct}{\cite}
\newcommand{\pr}{\protect\ref}
\newcommand{\su}{\subseteq}
\newcommand{\Dn}{\Delta^2_n}
\newcommand{\pa}{\partial}
\newcommand{\R}{{\Bbb R}}
\newcommand{\Z}{{\Bbb Z}}
\newcommand{\s}{{{\mathfrak{s}}^k}}
\newcommand{\p}{{\mathfrak{p}}}
\newcommand{\pp}{{\Bbb P_\p}}
\newcommand{\D}{\Delta}
\newcommand{\bs}{\bar{\sigma}}
\newcounter{numb}

\title{Topological embeddings into random 2-complexes}

\author{Michael Farber}
\thanks{Michael Farber was partially supported by a grant from the Leverhulme Foundation}
\address{School of Mathematical Sciences, 
	Queen Mary University of London,
	London E1 4NS, UK}
\email{M.Farber@@qmul.ac.uk}
\urladdr{\url{www.qmul.ac.uk/maths/profiles/farberm.html}}

\author{Tahl Nowik}
\address{Department of Mathematics, Bar-Ilan University, 
Ramat-Gan 5290002, Israel}
\email{tahl@@math.biu.ac.il}
\urladdr{\url{www.math.biu.ac.il/~tahl}}

\date{December 9, 2019}
\begin{abstract}
We consider 2-dimensional random simplicial complexes $Y$ in the multi-parameter model. We establish
the multi-parameter threshold for the property that every 2-dimensional simplicial 
complex $S$ admits a topological embedding into $Y$ asymptotically almost surely.
Namely, if in the procedure of the multi-parameter model, each $i$-dimensional simplex is taken independently  with probability $p_i=p_i(n)$, from a set of $n$ vertices, then the  threshold
is $p_0 p_1^3 p_2^2 =  \frac{1}{n}$.  
This threshold happens to coincide with the previously established thresholds for uniform hyperbolicity and triviality of the fundamental group. 

Our claim in one direction is in fact slightly stronger, namely, we show that if $p_0 p_1^3 p_2^2$ is sufficiently larger than $\frac{1}{n}$ 
then every $S$ has a fixed subdivision $S'$ which admits a simplicial embedding into $Y$ asymptotically almost surely.
The main geometric result we prove 
to this end is that given $\epsilon>0$, there is a subdivision $S'$ of $S$ such that 
every subcomplex $T \su S'$ has $\frac{f_0(T)}{f_1(T)}>\frac{1}{3}-\epsilon$ and $\frac{f_0(T)}{f_2(T)}>\frac{1}{2}-\epsilon$, 
where $f_i(T)$ denotes the number of simplices in $T$ of dimension $i$.

In the other direction we show that if $p_0 p_1^3 p_2^2$ is sufficiently smaller than $\frac{1}{n}$, 
then asymptotically almost surely, the torus does not admit a topological embedding into $Y$.
Here we use a result of Z. Gao which bounds the number of different triangulations of a surface.
\end{abstract}

\maketitle

\section{Introduction}\label{int}

The theory of random graphs is widely used for modeling large systems in various scientific applications. 
It has been recently realized that the higher dimensional generalization, namely, large random simplicial complexes,
offers even greater opportunities for mathematical modeling as it allows taking into account multiple interactions and not only pairwise relationships. The interest in \emph{large} complexes is manifested by studying the \emph{asymptotic} behavior of the random complexes, as the number of vertices grows.

The first model of random simplicial complexes was studied by Linial and Meshulam \cite{lm} and Meshulam and Wallach \cite{mw}. 
This model for simplicial complexes of dimension $r$ assumes that for a given vertex set, all simplices of dimension less than $r$ are present, and the top dimensional simplices are added independently with a given probability $p$. For the study of asymptotic properties, the probability $p$ is assumed to depend on the number $n$ of vertices. Many authors have contributed to the development of the Linial-Meshulam model by computing thresholds for various topological properties such as vanishing of the Betti numbers, simple connectivity, asphericity, and cohomological dimension of the fundamental group. 

In \cite{cf1}--\cite{cf3} Costa and Farber developed a more general model of random $r$-dimensional simplicial complexes. In this model a random simplicial complex depends on a string of probability parameters $(p_0, p_1, \dots, p_r)$, one for each dimension. 
In this multi-parameter setting, the thresholds for topological properties are not numbers depending on $n$, but rather boundaries of domains depending on $n$. In \cite{cf1}--\cite{cf3} domains corresponding to different asymptotic properties of the fundamental group are described. Another major result of \cite{cf1}--\cite{cf3} is the notion of a \emph{critical dimension} 
which fully portrays the properties of the Betti numbers in the multi-parameter model. 
The Linial-Meshulam model can be considered as a special case of the multi-parameter model by taking 
$p_0=p_1=\cdots=p_{r-1}=1$.

In \cite{fmn}, the authors studied an even more general model of random simplicial complexes, which is of two types, the \emph{lower} and \emph{upper} models, as described in Section \pr{bg}.  
The lower and upper models are Spanier-Whitehead dual to each other,
and the multi-parameter model mentioned above is a special case of the lower model. 
The notion of critical dimension for the upper model is also developed in \cite{fmn}.

In the present work we consider random 2-dimensional simplicial complexes $Y$ in the multi-parameter lower model. We find the threshold for the following property: 
\emph{Every 2-dimensional simplicial complex $S$ asymptotically almost surely admits a topological embedding into $Y$.}
We prove that the threshold for this property is $p_0 p_1^3 p_2^2 = \frac{1}{n}$, that is, if $p_0 p_1^3 p_2^2$ is sufficiently larger than $\frac{1}{n}$, namely, if there is an $\epsilon>0$ such that $np_0p_1^{3+\epsilon}p_2^{2+\epsilon}\to \infty$, then this property holds, whereas if $p_0 p_1^3 p_2^2$ is sufficiently smaller than $\frac{1}{n}$, namely, if $np_0p_1^{3}p_2^{2}\to 0$, then the property does not hold. 

The first half of this theorem generalizes a result of \cite{cc} regarding the Linial-Meshulam model, where it is proved for $p_2\geq n^{-\frac{1}{2}+\epsilon}$. This was strengthened by Gundert and Wagner in \cite{gw}, showing it for $p_2\geq cn^{-\frac{1}{2}}$. 
Versions of this result for higher dimensional complexes are also developed in \cite{gw}.
We also mention that our above threshold also happens to be the threshold for both uniform hyperbolicity and triviality of the fundamental group, as shown in \cite{cf2}.

\section{Background}\label{bg}

In this section we describe our model for random 2-dimensional simplicial complexes.
Let $\Dn$ denote the full 2-dimensional simplicial complex with vertex set $[n]=\{1,\dots,n\}$,
i.e.\ $\Dn$ includes all subsets of $[n]$ of size 1,2,3. They are regarded  as simplices of dimension 0,1,2, respectively.
A \emph{probability triple} $\p$ is a triple of sequences $(p_0(n),p_1(n),p_2(n))$ with $0 \leq p_i(n)\leq 1$.
We will however always suppress the argument $n$ from the expression $p_i(n)$ and simply write $p_i$.
A given argument $n$ will always be implied by the context.
A natural example for a probability triple is obtained as follows. Given a fixed $\alpha=(\alpha_0,\alpha_1,\alpha_2)$ with $\alpha_i\geq 0$, let $\p_\alpha$ be the probability triple given by 
$\p_\alpha=(n^{-\alpha_0},n^{-\alpha_1},n^{-\alpha_2})$.

We use a probability triple $\p$ to define a sequence $(\Omega_n,\Bbb P_\p)$ of probability spaces,
where $\Omega_n$ is the set of all simplicial subcomplexes of $\Dn$,
and the probability measure $\pp$ on $\Omega_n$ is defined as follows.
As an intermediate step we randomly select a subset $X \su \Dn$ which is not necessarily a simplicial complex.
Each $i$-dimensional simplex $\sigma \in \Dn$ is included in $X$ independently with probability $p_i$,
i.e.\ the probability for obtaining $X$ is 
$\prod_{\sigma \in X} p_{\dim\sigma} \prod_{\sigma \not\in X} q_{\dim\sigma}$, where $q_i=1-p_i$.
Now, the simplicial complex produced by this random process is the maximal simplicial complex $Y$ contained in $X$.
It is shown in \cite{cf1} and \cite{fmn} that the probability for obtaining any given simplicial complex $Y\su\Dn$ by this random process is
$$\pp(Y)=\prod_{\sigma \in Y} p_{\dim\sigma} 
\prod_{\sigma \in E(Y)} q_{\dim\sigma}$$
where $E(Y)=\{ \sigma\in\Dn : \sigma \not\in Y  \ \hbox{but} \   \pa \sigma \su Y\}$. 
Indeed, $Y$ is the maximal simplicial complex contained in $X$ iff every
$\sigma\in Y$ is in $X$ and every $\sigma\in E(Y)$ is not in $X$.
To see this, 
note that if $\sigma\in E(Y)$, then $Y \cup \{\sigma\}$ is also a simplicial complex. Also
note that if $\sigma \not\in Y$, then a minimal $\tau\su\sigma$ such that $\tau\not\in Y$,  is in $E(Y)$, (including when $\tau$ is a vertex). 
When we write $Y\in (\Omega_n, \pp)$ we will mean that $Y$ is randomly drawn from $\Omega_n$ using the probability measure $\pp$.

In \cite{fmn} this model is termed the \emph{lower model} for random simplicial complexes, in contrast to the \emph{upper model}, which looks at the \emph{minimal} simplicial complex \emph{containing} $X$. The setting in \cite{fmn} is more general, as the dimension of the simplices $\sigma$ is not restricted and each simplex $\sigma$ is assigned its own probability $p_\sigma$. 
 
For a 2-complex $S$ let $f_i(S)$, $i=0,1,2$, denote the number of simplices in $S$ of dimension $i$
(for us simplicial complexes are \emph{finite} by definition).
In this work we  will apply Theorem 1(B) of \ct{cf2} which states the following.
If $\p$ satisfies that for every subcomplex $T$ of $S$ 
$$\lim_{n \to\infty} n^{f_0(T)}p_0^{f_0(T)}p_1^{f_1(T)}p_2^{f_2(T)}=\infty,$$
then $S$ is \emph{simplicially} emebeddable into $Y\in (\Omega_n, \pp)$ a.a.s.\ (asymptotically almost surely, i.e.\ with probability converging to 1 as $n\to\infty$). 
In \ct{cf2} the corresponding claim for complexes of general dimension $r$ is proved.

\section{Statement of result}\label{sr}

In this work we look at random 2-complexes $Y \in (\Omega_n,\pp)$ for a given probability triple $\p=(p_0,p_1,p_2)$, and ask whether \emph{every} 2-complex $S$ is a.a.s.\ \emph{topologically} embeddable into $Y$.
We show that probability triples for which  $p_0 p_1^3 p_2^2 = \frac{1}{n}$ are the threshold for this property. That is,  triples $\p$ for which $p_0 p_1^3 p_2^2$ 
 is sufficiently larger than $\frac{1}{n}$ have the 
property that every 2-complex $S$ is a.a.s.\ topologically embeddable into 
$Y \in (\Omega_n,\pp)$, and triples $\p$ for which $p_0 p_1^3 p_2^2$ is sufficiently smaller than $\frac{1}{n}$ do not have this property. The precise statement of this claim appears in Theorem \pr{t1} below. It follows that for probability triples of the form $\p_\alpha$ 
the plane $\alpha_0 + 3\alpha_1 + 2\alpha_2=1$ is the threshold for this property. This is stated in
Corollary \pr{t2}.

\begin{thm}\label{t1}
Let  $\p=(p_0,p_1,p_2)$ be a probability triple.
\begin{enumerate}
\item 
 If 
$$\lim_{n\to\infty}n p_0 p_1^{3+\epsilon}p_2^{2+\epsilon}=\infty$$ for some fixed $\epsilon >0$,
then every 2-dimensional simplicial complex $S$ has a simplicial subdivision $S'$ which is a.a.s.\  simplicially embeddable into $Y\in(\Omega_n,\pp)$. In particular, $S$ is a.a.s.\ toplogically  embeddable into $Y\in(\Omega_n,\pp)$. 
\item
If $$\lim_{n\to\infty}n p_0 p_1^3p_2^2=0$$  then the torus is a.a.s.\ not topologically embeddable into $Y\in(\Omega_n,\pp)$.
\end{enumerate}
\end{thm}

\begin{cor}\label{t2} 
Let  $\alpha=(\alpha_0,\alpha_1,\alpha_2)$ with $\alpha_i\geq 0$ 
and let $\p_\alpha=(n^{-\alpha_0},n^{-\alpha_1},n^{-\alpha_2})$.
\begin{enumerate}
\item  If $\alpha_0 + 3\alpha_1 + 2\alpha_2<1$ 
then every 2-dimensional simplicial complex $S$ has a 
simplicial subdivision $S'$ which is a.a.s.\ simplicially embeddable into 
$Y\in (\Omega_n, \Bbb P_{\p_\alpha})$.
In particular, $S$ is a.a.s.\ toplogically  embeddable into $Y\in (\Omega_n, \Bbb P_{\p_\alpha})$.
\item If $\alpha_0 + 3\alpha_1 + 2\alpha_2>1$ then the torus is a.a.s.\ not topologically embeddable into $Y\in (\Omega_n, \Bbb P_{\p_\alpha})$.
\end{enumerate}
\end{cor}

The main effort of this work is the proof of Proposition \pr{m} below, for which we use a particular subdivision scheme for 2-complexes which is suitable for our purposes. The first subdivision is barycentric, but all subsequent subdivisions divide each triangle into four triangles as shown in Figure \pr{f1}. Geometrically, we view the $k$th subdivision $\s \D$ of a 2-simplex $\D$ as a regular hexagon divided into $6\cdot 4^{k-1}$ equilateral triangles.

\begin{figure}
\begin{tikzpicture}[scale=0.6]

\draw(0,0)++(60:2)--(0:4) 
(0,0) ++(-60:2) --(0:4)
(60:2)--(-60:2)
;
\filldraw (0,0) ++(60:2) circle(3pt) 
(0,0)++(0:4) circle(3pt)
(0,0) ++(-60:2) circle(3pt) 
;

\draw[->][line width=2] (4.5,0)--(11/2,0);

\begin{scope}[shift={(5,0)}]
\draw(0,0)++(60:2)--(0:4) 
(0,0) ++(-60:2) --(0:4)
(60:2)--(-60:2)
;
\draw (1,0)--(4,0)
(0,0)++(60:2)--++(-60:3)
(0,0)++(-60:2)--++(60:3)

;

\filldraw (0,0) ++(60:2) circle(3pt) 
(0,0)++(0:4) circle(3pt) node[right] {$=$} 
(0,0) ++(-60:2) circle(3pt);

\end{scope}

\begin{scope}[shift={(10,0)}]
\draw (0,0) -- ++(0:4)
(0,0)++(60:2)--++(0:2)
(0,0)++(-60:2)--++(0:2)

(0,0)--++(60:2)
(0,0)++(-60:2)--++(60:4)
(0,0)++(-60:2)++(0:2)--++(60:2)

(0,0)--++(-60:2)
(0,0)++(60:2)--++(-60:4)
(0,0)++(60:2)++(0:2)--++(-60:2)
;
\filldraw (0,0) ++(60:2) circle(3pt) 
(0,0)++(0:4) circle(3pt)
(0,0) ++(-60:2) circle(3pt) 
;
\draw[->][line width=2] (4.5,0)--(11/2,0);
\end{scope}

\begin{scope}[shift={(16,0)}]
\draw (0,0) -- ++(0:4)
(0,0)++(60:1)--++(0:3)
(0,0)++(60:2)--++(0:2)
(0,0)++(-60:1)--++(0:3)
(0,0)++(-60:2)--++(0:2)
(0,0)--++(60:2)
(0,0)++(-60:1)--++(60:3)
(0,0)++(-60:2)--++(60:4)
(0,0)++(-60:2)++(0:1)--++(60:3)
(0,0)++(-60:2)++(0:2)--++(60:2)
(0,0)--++(-60:2)
(0,0)++(60:1)--++(-60:3)
(0,0)++(60:2)--++(-60:4)
(0,0)++(60:2)++(0:1)--++(-60:3)
(0,0)++(60:2)++(0:2)--++(-60:2)
;
\filldraw (0,0) ++(60:2) circle(3pt) 
(0,0)++(0:4) circle(3pt)
(0,0) ++(-60:2) circle(3pt) 
;
\draw[->][line width=2] (4.5,0)--(11/2,0);
\end{scope}

\begin{scope}[shift={(22,0)}]
\draw (0,0) -- ++(0:8/2)
(0,0)++(60:1/2)--++(0:7/2)
(0,0)++(60:2/2)--++(0:6/2)
(0,0)++(60:3/2)--++(0:5/2)
(0,0)++(60:4/2)--++(0:4/2)
(0,0)++(-60:1/2)--++(0:7/2)
(0,0)++(-60:2/2)--++(0:6/2)
(0,0)++(-60:3/2)--++(0:5/2)
(0,0)++(-60:4/2)--++(0:4/2)

(0,0)--++(60:4/2)
(0,0)++(-60:1/2)--++(60:5/2)
(0,0)++(-60:2/2)--++(60:6/2)
(0,0)++(-60:3/2)--++(60:7/2)
(0,0)++(-60:4/2)--++(60:8/2)
(0,0)++(-60:4/2)++(0:1/2)--++(60:7/2)
(0,0)++(-60:4/2)++(0:2/2)--++(60:6/2)
(0,0)++(-60:4/2)++(0:3/2)--++(60:5/2)
(0,0)++(-60:4/2)++(0:4/2)--++(60:4/2)

(0,0)--++(-60:4/2)
(0,0)++(60:1/2)--++(-60:5/2)
(0,0)++(60:2/2)--++(-60:6/2)
(0,0)++(60:3/2)--++(-60:7/2)
(0,0)++(60:4/2)--++(-60:8/2)
(0,0)++(60:4/2)++(0:1/2)--++(-60:7/2)
(0,0)++(60:4/2)++(0:2/2)--++(-60:6/2)
(0,0)++(60:4/2)++(0:3/2)--++(-60:5/2)
(0,0)++(60:4/2)++(0:4/2)--++(-60:4/2)
;
\filldraw (0,0) ++(60:2) circle(3pt) 
(0,0)++(0:4) circle(3pt)
(0,0) ++(-60:2) circle(3pt) 
;
\end{scope}
\end{tikzpicture}
\caption{$\D$ \  $\longrightarrow$ \ two depictions of ${\mathfrak{s}}^1 \D $ \ $\longrightarrow$ \ ${\mathfrak{s}}^2 \D$ \ $\longrightarrow$ \  ${\mathfrak{s}}^3 \D$}\label{f1}
\end{figure}
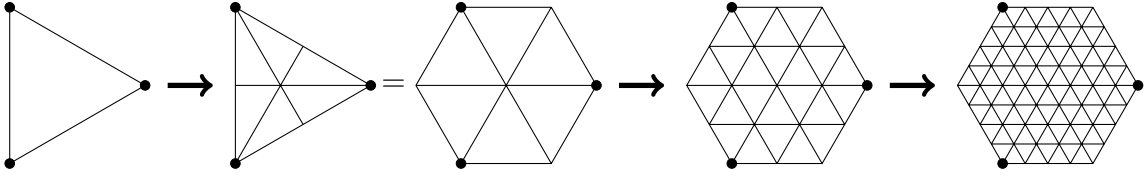

For a 2-complex $S$ let $f_i(S)$, $i=0,1,2$, denote the number of simplices in $S$ of dimension $i$, and define $\mu_i(S)=\frac{f_0(S)}{f_i(S)}$, $i=1,2$, whenever $f_i(S) \neq 0$.

\begin{prop}\label{m}
Let $S$ be a 2-complex, and let $\epsilon>0$. Then for sufficiently large $k$, 
every subcomplex $T$ of $\s S$ has
\begin{enumerate}
\item $\mu_1(T)>\frac{1}{3}-\epsilon$.
\item $\mu_2(T)>\frac{1}{2}-\epsilon$.
\end{enumerate}
\end{prop}

The proofs of parts (1) and (2) of Proposition \pr{m} appear in Sections \pr{ed} and \pr{tr} respectively.
Once Proposition \pr{m} is established, the proof of Theorem \pr{t1}(1) is immediate. Indeed,
by Proposition \pr{m}, given $\epsilon >0$ there exists a subdivision $S'$ of $S$ such that every subcomplex $T$ of $S'$ satisfies 
$\mu_1(T)>\frac{1}{3+\epsilon}$ and 
$\mu_2(T)>\frac{1}{2+\epsilon}$. For every subcomplex $T$ of $S'$ we have
$$n^{f_0(T)}p_0^{f_0(T)}p_1^{f_1(T)}p_2^{f_2(T)}=
\bigg(n p_0 p_1^{\frac{1}{\mu_1(T)}}p_2^{\frac{1}{\mu_2(T)}}\bigg)^{f_0(T)}\geq 
\bigg(n p_0 p_1^{3+\epsilon}p_2^{2+\epsilon}\bigg)^{f_0(T)}\ \mathop{\longrightarrow}_{n\to\infty}\ \infty.$$
By Theorem 1(B) of \ct{cf2} stated in the previous section, $S'$ is a.a.s.\ simplically embeddable into $Y\in(\Omega,\pp)$, completing the proof of Theorem \pr{t1}(1). The proof of Theorem \pr{t1}(2) appears in Section \pr{ne}. Corollary \pr{t2} easily follows from Theorem \pr{t1}.

\section{Domains in the open 2-simplex}\label{dm}

In this section we focus our attention on one 2-simplex $\D$ of the original complex $S$. 
We will be interested in open planar surfaces of the form $\D-(\pa\D\cup T)$ where $T$ is a pure subcomplex of $\s\D$. By a pure subcomplex we mean a subcomplex which is the union of simplices of $\s\D$ of fixed dimension (either 1 or 2).  

\begin{dfn}
Given a simplicial complex $S$ and a subcomplex $T$ of $S$, the \emph{open complex defined by $S-T$} is the collection $W$ of all \emph{open} simplices in $S-T$. The underlying space of $W$, i.e.\ the union of open simplices in $W$, will be denoted $|W|$.
\end{dfn}

\begin{dfn}
Let $\D$ be a fixed  2-simplex.
\begin{enumerate}
\item We denote by $V_k$ the open complex defined by $\s\D -\pa(\s\D)$, 
where $\pa(\s\D)$ is the 
%1-dimensional 
subcomplex of  $\s\D$ which is the topological boundary of $\D$.
So, $|V_k|$ is the interior of $\D$, i.e.\ an open disc. The open simplices of dimensions 0,1,2 in $V_k$ will be called vertices, edges and triangles, respectively. 
\item Given a collection $L=\{E_1,\dots,E_m\}$ of \emph{closed} simplices in $\s\D$ of fixed dimension, we denote by $U_L$
the open complex defined by
$$\s\D -\bigg(\pa(\s\D)\cup\bigcup_{E_i\in L}E_i\bigg).$$ 
An open complex of the form $U_L$ will be called a \emph{domain} in $V_k$. 
The fixed dimension of the simplices in $L$ will be called the ``type'' of $U_L$.
Note that $V_k$ itself is a domain with $L=\varnothing$. 
\item For a domain $U=U_L$ we define $f_i(U)$
to be the number of open simplices of dimension $i$ in $U$, $i=0,1,2$, and define $\mu_i(U)=\frac{f_0(U)}{f_i(U)}$, $i=1,2$, whenever $f_i(U)\neq 0$. 
\end{enumerate}
\end{dfn}

\begin{lemma}\label{01}
Let $U=U_L$ be a domain in $V_k$, then $\mu_1(U)<\frac{1}{3}$ and  
$\mu_2(U)<\frac{1}{2}$ whenever they are defined.
\end{lemma}

\begin{pf}
If there are no vertices in $U$ then $\mu_1(U)=\mu_2(U)=0$ and we are done. Otherwise, let $v$ be a vertex in $U$.
Since $U$ is open, all 6 edges surrounding $v$ must be in $U$. On the other hand, 
given an edge $e$ in $U$, at most 2 vertices in $U$ lie at its endpoints. This gives $\mu_1(U)\leq\frac{1}{3}$, and to show strict inequality we must show that
if $\mu_1(U)$ is defined, i.e.\ if there are edges in $U$, then at least one of them does not have both endpoints in $U$. Indeed, starting with an edge $e$ in $U$, follow a path of edges in $V_k$ to $\pa V_k$. 
Since $U$ is open, this path must include an edge with the desired property.

Similarly, if $v$ is a vertex in $U$ then since $U$ is open, all 6 triangles surrounding $v$ are in $U$. On the other hand at most 3 vertices in $U$ lie on the boundary of a given triangle in $U$, and triangles adjacent to the boundary of $U$ have strictly less than 3 such vertices. This gives  $\mu_2(U)<\frac{1}{2}$.
\end{pf}

Regarding $V_k$ itself we note the following.

\begin{prop}\label{vk}
$\lim_{k\to\infty} \mu_1(V_k)=\frac{1}{3}$ \ and \ $\lim_{k\to\infty} \mu_2(V_k)=\frac{1}{2}$.
\end{prop}

\begin{pf} 
We claim that for $k \geq 1$,
 $f_0(V_k)= 3\cdot 4^{k-1} -3\cdot 2^{k-1}   +1   $,
$f_1(V_k)= 9\cdot 4^{k-1} -3\cdot 2^{k-1} $ and 
$f_2(V_k)=6\cdot 4^{k-1}$, from which the statement follows.
This is verified by induction noting that
$f_0(V_1)=1$, $f_1(V_1)=6$, $f_2(V_1)=6$, and for $k\geq 1$,
$f_0(V_{k+1})=f_0(V_k)+f_1(V_k)$, 
$f_1(V_{k+1})=2f_1(V_k)+3f_2(V_k)$ and
$f_2(V_{k+1})=4f_2(V_k)$.
\end{pf}

An open hexagon in $V_k$ is a domain in $V_k$ which is convex and has six sides, which are necessarily parallel to the six sides of $V_k$. Two examples appear in Figure \pr{f4}. We denote the boundary of an open hexagon $U$ by $\pa U$, though strictly speaking it is the boundary of the corresponding closed hexagon.

\begin{lemma}\label{hx}
For $U$ an open hexagon in $V_k$, let $f_i=f_i(U)$ and $b=$ the number of edges in $\pa U$, then 
 $f_0=\frac{f_2-b+2}{2}$
and $f_1=\frac{3f_2-b}{2}$.
\end{lemma}

\begin{pf}
Every triangle in $U$ is surrounded by three edges. In this way each edge in $U$ is counted twice, and each edge in $\pa U$ is counted once, so $3f_2=2f_1+b$. Also, since the boundary of the closed hexagon contributes 0 to its Euler characteristic, we have $f_0-f_1+f_2=1$. Our claim follows by combining these two equalities.
\end{pf}

For a domain $U$ in $V_k$, a ``shifted image'' of $U$ in $V_k$ is a domain in $V_k$ obtained from $U$ by a (perhaps trivial) translation. A ``side'' of an open hexagon $U$ in $V_k$ is one of the six sides of $\pa U$.

\begin{lemma}\label{hx2}
Let $U$ be an open hexagon in $V_k$ with $U\neq V_k$. Then there is a shifted image of $U$ in $V_k$  
whose longest side is not included in $\pa V_k$. 
\end{lemma}

\begin{pf}
If $I$ is one of the six sides of $U$, then we define its length $\ell(I)$ to be the number of edges in $I$.
Let $I$ be the longest side of $U$, let $t=\ell(I)$, and let $s$ denote the side length of $V_k$.
If $t\geq s+1$ we are done since $I$ cannot be included in $\pa V_k$.
If $t\leq s-1$ then $U$ is included in a regular hexagon of side $s-1$ and so the claim is again clear.
So assume $t=s$.
It cannot be that all sides of $U$ are of length $s$ since then $U=V_k$. So the sides of $U$ may be consecutively labeled $I_1,\dots,I_6$ so that $\ell(I_1)=s$, $\ell(I_2)\leq s-1$, and $\ell(I_3)\leq s$. 
If $I_1$ is not included in $\pa V_k$ we are done, so assume $I_1$ is a side of $V_k$. 
Now let $U'$ be the hexagon with corresponding sides $J_1,\dots,J_6$ with $J_1=I_1$ so $\ell(J_1)=s$, and having $\ell(J_2)=s-1$ and $\ell(J_3)=\ell(J_4)=s$, as seen in Figure \pr{f2}. It is easy to see that necessarily 
$U\su U'$ and so $U$ may be shifted by one to the right, pulling $I_1$ away from $\pa V_k$.
\end{pf}

\begin{figure}
\begin{center}
\begin{tikzpicture}[scale=0.8]
\fill[gray!35]
(0,0)--++(-60:4)--++(0:3)--++(60:4)--++(120:4)--++(180:3)
; 
\draw (0,0) -- ++(0:8)
(60:1)--++(0:7)
(60:2)--++(0:6)
(60:3)--++(0:5)
(60:4)--++(0:4)
(-60:1)--++(0:7)
(-60:2)--++(0:6)
(-60:3)--++(0:5)
(-60:4)--++(0:4)

(0,0)--++(60:4)
(-60:1)--++(60:5)
(-60:2)--++(60:6)
(-60:3)--++(60:7)
(-60:4)--++(60:8)
(-60:4)++(0:1)--++(60:7)
(-60:4)++(0:2)--++(60:6)
(-60:4)++(0:3)--++(60:5)
(-60:4)++(0:4)--++(60:4)

(0,0)--++(-60:4)
(60:1)--++(-60:5)
(60:2)--++(-60:6)
(60:3)--++(-60:7)
(60:4)--++(-60:8)
(60:4)++(0:1)--++(-60:7)
(60:4)++(0:2)--++(-60:6)
(60:4)++(0:3)--++(-60:5)
(60:4)++(0:4)--++(-60:4)

(-60:2)++(180:1) node {$I_1=J_1$}
(-60:4)++(0:1.7)++(-90:0.4) node {$J_2$}
(0:7.3)++(-120:1.7) node {$J_3$}
(0:7.3)++(120:2.7) node {$J_4$}
;

\end{tikzpicture}
\end{center}
\caption{The hexagon $U'$  with $\ell(J_2)=s-1$ and $\ell(J_1)=\ell(J_3)=\ell(J_4)=s$.}\label{f2}
\end{figure}
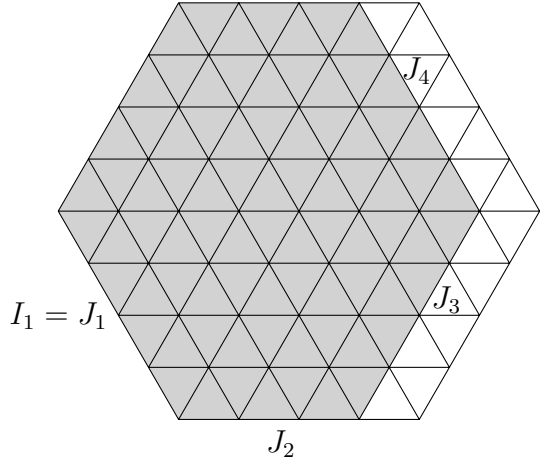

We conclude this section with the following trivial lemma that we will use repeatedly throughout this work. 
It will be used e.g.\ when say we have disjoint $A,B$ with $\mu_i(A)>c$ and $\mu_i(B)>c$ to deduce $\mu_i(A\cup B)>c$.
Another instance where it will be used is when say $f_0=x$, $f_i=y$, 
and we are interested to know whether  $\mu_i$ increases if we replace $x,y$ by 
$x\pm a, y\pm b$. 

\begin{lemma}\label{p} \ 
\begin{enumerate}
\item Let $x,a\geq 0$,  $y,b>0$. Then $\frac{x+a}{y+b}=t\frac{x}{y}+(1-t)\frac{a}{b}$ with $0<t<1$. 
That is, either $\frac{x}{y}=\frac{a}{b}=\frac{x+a}{y+b}$, or
$\frac{x+a}{y+b}$ is strictly between $\frac{x}{y}$ and $\frac{a}{b}$.
In particular if $\frac{a}{b}>\frac{x}{y}$ then $\frac{x+a}{y+b}>\frac{x}{y}$.
\item Let $x\geq a\geq 0$,  $y>b>0$. Then $\frac{x}{y}=t\frac{x-a}{y-b}+(1-t)\frac{a}{b}$ with $0<t<1$. That is,  either $\frac{x}{y}=\frac{a}{b}=\frac{x-a}{y-b}$, or
$\frac{x}{y}$ is strictly between $\frac{x-a}{y-b}$ and $\frac{a}{b}$.
In particular if $\frac{a}{b} < \frac{x}{y}$ then $\frac{x-a}{y-b}>\frac{x}{y}$.
\end{enumerate}
\end{lemma}

\begin{pf} 
For (1) take $t=\frac{y}{y+b}$. For (2) take $t=\frac{y-b}{y}$.
\end{pf}

\section{Counting triangles}\label{tr}

In this section we prove Proposition \pr{m}(2), where
the main technical effort is the following.

\begin{prop}\label{u}
For sufficiently large $k$, every nonempty type-2 domain $U=U_L$  in $V_k$ with $U \neq V_k$ 
has $\mu_2(U) < \mu_2(V_k)$.
\end{prop}

\begin{pf} 
Throughout this proof, a type-2 domain will simply be called a domain.
We establish our claim by showing that for every domain $U \neq V_k$ there is another domain $U'$ with $\mu_2(U')>\mu_2(U)$.
If $\mu_2(U)=0$ then simply take $U'=V_k$, so we assume $\mu_2(U)>0$, i.e.\ $f_0(U)>0$, which implies $f_2(U)\geq 6$.

To each vertex $v$ of $V_k$ we assign a cyclically ordered 6-tuple of symbols $x$ and $o$ describing the six triangles surrounding $v$, where an $x$ denotes a 2-simplex in $L$ and an $o$ denotes a 2-simplex which is not in $L$, and so its interior (and perhaps more of it) is included in $U$.  
To each vertex of $\pa V_k$ we also attach such a 6-tuple. It will include $x$'s and $o$'s defined as before, and then completed with a sequence of $x$'s to reach a total of six symbols. One may think of these additional ``virtual'' $x$'s as corresponding to triangles outside $V_k$ in a tiling of a neighborhood of $V_k$ in the plane. See Figure \pr{f3}, where the triangles in $U$ are depicted in gray color.

\begin{figure}
\begin{center}
\begin{tikzpicture}
\fill[gray!35]
(-60:1)--++(-60:1)--++(0:1)--++(60:1)
(0,0)--++(-60:1)--++(60:1)--++(0:1)--++(120:1)--++(180:1)
(0:2)--++(0:2)--++(120:1)--++(180:1)
(60:1)++(0:1)--++(0:1)--++(120:1)
; 

\draw (0,0) -- ++(0:4)
(60:1)--++(0:3)
(60:2)--++(0:2)
(-60:1)--++(0:3)
(-60:2)--++(0:2)
(0,0)--++(60:2)
(-60:1)--++(60:3)
(-60:2)--++(60:4)
(-60:2)++(0:1)--++(60:3)
(-60:2)++(0:2)--++(60:2)
(0,0)--++(-60:2)
(60:1)--++(-60:3)
(60:2)--++(-60:4)
(60:2)++(0:1)--++(-60:3)
(60:2)++(0:2)--++(-60:2)
;
\draw (-60:1)--++(-120:0.5)
(-60:1)--++(-180:0.5)

(0:4)--++(60:0.5)
(0:4)--++(0:0.5)
(0:4)--++(-60:0.5)

(-60:1)++(-90:0.4) node {$x$}
(-60:1)++(-150:0.4) node {$x$}
(-60:1)++(150:0.4) node {$x$}
(-60:1)++(-30:0.4) node {$o$}
(-60:1)++(30:0.4) node {$x$}
(-60:1)++(90:0.4) node {$o$}

(60:1)++(0:2)++(-90:0.4) node {$o$}
(60:1)++(0:2)++(-150:0.4) node {$x$}
(60:1)++(0:2)++(150:0.4) node {$o$}
(60:1)++(0:2)++(-30:0.4) node {$o$}
(60:1)++(0:2)++(30:0.4) node {$x$}
(60:1)++(0:2)++(90:0.4) node {$x$}

(0:4)++(-90:0.4) node {$x$}
(0:4)++(-150:0.4) node {$x$}
(0:4)++(150:0.4) node {$o$}
(0:4)++(-30:0.4) node {$x$}
(0:4)++(30:0.4) node {$x$}
(0:4)++(90:0.4) node {$x$}

;
\filldraw 
(-60:1) circle(2pt) 
(60:1)++(0:2) circle(2pt)
(0:4) circle(2pt)
;

\end{tikzpicture}
\end{center}
\caption{Cyclically ordered 6-tuples of $o$'s and $x$'s for type-2 domain.}\label{f3}
\end{figure}
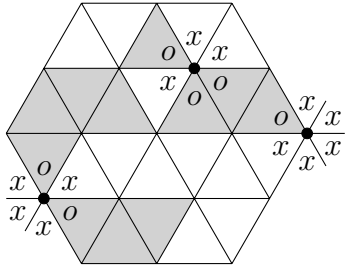

If for some vertex $v$ of $V_k$ the corresponding 6-tuple includes one or two $x$'s
then we may remove the corresponding 2-simplices from $L$.
This adds one or two triangles to $U$, and at least one vertex, namely $v$ itself. 
So if $a$ is the number of additional vertices and $b$ is the number of additional triangles then $\frac{a}{b}\geq\frac{1}{2}$. From Lemmas \pr{01} and \pr{p}(1) it then follows that $\mu_2$ strictly increases. So we may assume that every 6-tuple includes either 0 or at least 3 $x$'s. This is also true for the vertices of $\pa V_k$ since each such vertex has 3 or 4 virtual $x$'s, as seen in Figure \pr{f3}.

If for some vertex $v$ of $V_k$ or $\pa V_k$ the corresponding 6-tuple includes a fragment of the form $xox$ then we may add the middle 2-simplex to $L$, changing the given fragment to $xxx$. This removes one triangle from $U$, and leaves the number of vertices unchanged, and so it strictly increases $\mu_2$. 

We thus need to establish our claim when all 6-tuples have 0 or at least 3 $x$'s, and a fragment of the form $xox$ does not appear. It follows that all 6-tuples that include both $x$'s and $o$'s are either  $ooxxxx$ or $oooxxx$. That is, all boundary angles of the connected components of $U$ are either $120^\circ$ or $180^\circ$. It follows that each connected component of $U$ is an open hexagon. 
Indeed, note that a connected component of $U$ cannot have more than one boundary component, since then it would have an interior boundary component, along which there would be angles greater than $180^\circ$.
By Lemma \pr{p}(1) it is enough to prove our claim for each such hexagon separately, and so we assume $U$ is a single hexagon.

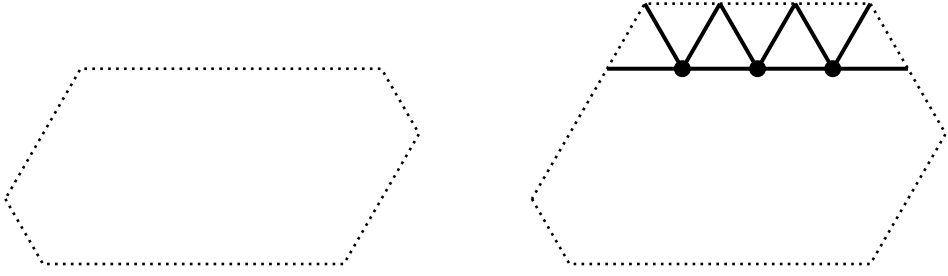
\begin{figure}
\begin{center}
\begin{tikzpicture}

\draw[line width=1,dotted] (0,0)--++(-60:1)--++(0:4)--++(60:2)--++(120:1)--++(180:4)--++(-120:2)
;

\begin{scope}[shift={(7,0)}]

\draw[line width=1,dotted] (0,0)--++(-60:1)--++(0:4)--++(60:2)--++(120:2)--++(180:3)--++(-120:3)
;
\draw[line width=1.5]
(60:2)--++(0:4)
(60:3)--++(-60:1)
(60:3)++(0:1)--++(-60:1)
(60:3)++(0:2)--++(-60:1)
(60:2)++(0:1)--++(60:1)
(60:2)++(0:2)--++(60:1)
(60:2)++(0:3)--++(60:1)
;

\filldraw 
(0,0)++(60:2)++(0:1) circle(3pt) 
(0,0)++(60:2)++(0:2) circle(3pt)
(0,0)++(60:2)++(0:3) circle(3pt)
;
\end{scope}
\end{tikzpicture}
\end{center}
\caption{Adding a row of triangles along the longest side of an open hexagon.}\label{f4}
\end{figure}

Let $I$ be the longest side of $U$, having length $t$. By Lemma \pr{hx2} we may shift $U$ in $V_k$ so that $I$ is not included in the boundary of $V_k$.
We may now look at a larger hexagon $U'$ in $V_k$ obtained from $U$ by adding a full row of triangles along $I$ as in Figure \pr{f4}. This process adds $t-1$ vertices and $2t-1$ triangles, as seen in Figure \pr{f4} where $t=4$.
By Lemma \pr{hx} we have  
$\mu_2(U)=\frac{f_2-b+2}{2f_2}$ where $b$ is the number of edges in $\pa U$.
We would like to use Lemma \pr{p}(1) to show that $\mu_2(U')>\mu_2(U)$, so we would like to have $\frac{t-1}{2t-1} > \frac{f_2-b+2}{2f_2}$, i.e.\ $f_2 < (b-2)(2t-1)$.

We show that this indeed 
holds if $b\geq 12$. Since $I$ is the longest side, we have $t\geq \frac{b}{6}$
and so it is enough to show that $f_2 < (b-2)(\frac{b}{3}-1)=\frac{1}{3}(b^2-5b+6)$. 
If we think of our triangles as equilateral triangles in the plane with sides of length 1, then the area $A$ of $U$ is $\frac{\sqrt{3}}{4}f_2$, the length of $\pa U$ is $b$,  and our required inequality may be rephrased as $A < \frac{1}{4\sqrt{3}}(b^2-5b+6)$. By the isoperimetric inequality in the plane we 
know $A \leq \frac{b^2}{4\pi}$ and so it is enough to have  $\frac{b^2}{4\pi} < \frac{1}{4\sqrt{3}}(b^2-5b)$, i.e.\ $\frac{\sqrt{3}}{\pi} < 1-\frac{5}{b}$, which holds for $b\geq 12$.

Finally, there are only finitely many possible hexagons with $b\leq 11$, and so by Lemmas \pr{01} and \pr{vk},
for sufficiently large $k$ we have $ \mu_2(V_k) > \mu_2(U)$ for every hexagon $U$ with $b\leq 11$, and so for such $k$ our statement is established.
\end{pf}

We may now proceed to prove Proposition \pr{m}(2), namely, that given $\epsilon>0$, then for sufficiently large $k$, every subcomplex $T$ of $\s S$ has $\mu_2(T)>\frac{1}{2}-\epsilon$.
(Recall that $T$ is a subcomplex in the usual sense, i.e.\ closed in $S$.) 
We first note that it is enough to prove the claim for pure 2-dimensional subcomplexes since given any subcomplex $T$ we may replace it by the subcomplex which is the union of its (closed) 2-simplices. This does not change $f_2(T)$ and may only reduce $f_0(T)$ and so $\mu_2(T)$ may only be reduced.
Now, given $\epsilon > 0$, by Proposition \pr{vk} we have   $\mu_2(V_k)>\frac{1}{2}-\epsilon$ for sufficiently large $k$, and we show that every pure 2-dimensional subcomplex $T$ of $\s S$ satisfies 
$\mu_2(T)>\frac{1}{2}-\epsilon$. 
For a fixed 2-simplex $\D$ of $S$ we look at the corresponding $V_k$. If $T \cap V_k \neq \varnothing$ we claim that 
$\mu_2(T \cap V_k) > \frac{1}{2}-\epsilon$ (where $\mu_2(T \cap V_k)$ is defined, as usual, by counting open simplices). Indeed this holds if $T \cap V_k = V_k$. Otherwise look at $U=V_k-(T\cap V_k)$, which is a nonempty type-2 domain in $V_k$. 
By Proposition \pr{u} we have $\mu_2(U) < \mu_2(V_k)$, and so by Lemma \pr{p}(2) we have 
$\mu_2(T \cap V_k) = \mu_2(V_k-U)> \mu_2(V_k) > \frac{1}{2}-\epsilon$. If we look at the union of  $T \cap V_k$ over
all 2-simplices of $S$ then by Lemma \pr{p}(1) we still have $\mu_2>\frac{1}{2}-\epsilon$. 
All that is missing from the full subcomplex $T$ is $T \cap \s S^1$, where $S^1$ is the 1-skeleton of $S$.
But $T \cap \s S^1$ may only contribute additional vertices, and no additional 2-simplices, so $\mu_2$ may only increase, and so we finally get 
$\mu_2(T)>\frac{1}{2}-\epsilon$, completing the proof of Proposition \pr{m}(2).

\subsection*{Discussion}
At this point the reader may have noticed that all our computations up to and including the proof of Proposition \pr{u} dealt with domains. Only in the previous paragraph have we finally returned to subcomplexes $T$ of $S$, and domains made their appearance as the \emph{complement} of $T\cap V_k$ in $V_k$.
The reason for concentrating all our attention on the complement of $T\cap V_k$ in $V_k$, rather than on $T$ itself, is the following.

Our aim was to show that $\mu_2(T)$ is sufficiently large. We have done this by separately considering the intersection of $T$ with each 2-simplex of $S$. If we had considered $T\cap \Delta$ for the \emph{closed} simplex $\Delta$, then the vertices in $T \cap \s S^1$ would have been multiply counted since each 
1-simplex of $S$ may be included in many 2-simplices of $S$. For this reason we were forced to consider
$T \cap int \Delta = T\cap V_k$. This is a peculiar geometric object with only part of its boundary present. We were not able however to ignore this partial boundary since without it we would be left with a domain in $V_k$, for which we have seen above that in fact $\mu_2$ tends to be \emph{small}.
On the other hand, the fact that the \emph{complement} of $T\cap V_k$ in $V_k$ is a domain, and thus has small $\mu_2$, was to our advantage since we could use it, via Lemma \pr{p}(2), to deduce that $T\cap V_k$ itself has large $\mu_2$.

The same approach will be used in the next section, for proving Proposition \pr{m}(1).

\section{Counting edges}\label{ed}
In this section we prove Proposition \pr{m}(1). Again, the main technical effort is the following.

\begin{prop}\label{u1}
For sufficiently large $k$, every type-1 domain $U=U_L$  in $V_k$ with $U \neq V_k$ 
has $\mu_1(U) < \mu_1(V_k)$, if $\mu_1(U)$ is defined.
\end{prop}

\begin{pf} 
Throughout this proof, a type-1 domain will simply be called a domain.
Again, we establish our claim by showing that for every domain $U \neq V_k$ there is another domain $U'$ with $\mu_1(U')>\mu_1(U)$.
If $\mu_1(U)=0$ then simply take $U'=V_k$, so we assume $\mu_1(U)>0$, i.e.\ $f_0(U)>0$, which implies $f_1(U)\geq 6$. Note that since the simplices in $L$ are of dimension 1, 
all triangles of $V_k$ are in $U$.

To each vertex $v$ of $V_k$ we assign a cyclically ordered 6-tuple of symbols $x$ and $o$,
this time describing the six \emph{edges} surrounding $v$, where an $o$ denotes an edge of $V_k$ that is in $U$, and an $x$ denotes an edge of $V_k$ which is not in $U$.  
For each vertex of $\pa V_k$ we also define such a 6-tuple, by adding, as before, a sequence of ``virtual'' $x$'s, to reach a total of six symbols. See Figure \pr{f5}, where the edges of $U$ are depicted in white color (i.e.\ they are missing).
Since all triangles of $V_k$ are in $U$, and they are also white, the regions appearing in Figure \pr{f5} are in fact the connected components of $U$. Note that a 6-tuple corresponding to a vertex of $\pa V_k$ includes at most 2 $o$'s, and so at least 4 $x$'s.

\begin{figure}
\begin{center}
\begin{tikzpicture}

\draw (0,0) -- ++(0:2)
(60:2)--++(0:2)
(-60:2)--++(0:2)
(0,0)--++(60:2)
(-60:1)--++(60:2)
(-60:2)--++(60:4)
(-60:2)++(0:1)++(60:2)--++(60:1)
(-60:2)++(0:2)--++(60:2)
(0,0)--++(-60:2)
(60:1)--++(-60:3)
(60:2)++(0:1)--++(-60:3)
(60:2)++(0:2)--++(-60:2)

(-60:1)--++(-120:0.5)
(-60:1)--++(-180:0.5)

(0:4)--++(60:0.5)
(0:4)--++(0:0.5)
(0:4)--++(-60:0.5)

;

 \draw
(-60:1)++(0:0.4) node {$o$}
(-60:1)++(60:0.4) node {$x$}
(-60:1)++(120:0.4) node {$x$}
(-60:1)++(180:0.4) node {$x$}
(-60:1)++(-60:0.4) node {$x$}
(-60:1)++(-120:0.4) node {$x$}

(60:1)++(0:2)++(0:0.4) node {$o$}
(60:1)++(0:2)++(60:0.4) node {$x$}
(60:1)++(0:2)++(120:0.4) node {$x$}
(60:1)++(0:2)++(180:0.4) node {$o$}
(60:1)++(0:2)++(-60:0.4) node {$x$}
(60:1)++(0:2)++(-120:0.4) node {$x$}

(0:4)++(0:0.4) node {$x$}
(0:4)++(60:0.4) node {$x$}
(0:4)++(120:0.4) node {$x$}
(0:4)++(180:0.4) node {$o$}
(0:4)++(-60:0.4) node {$x$}
(0:4)++(-120:0.4) node {$x$}

;

\filldraw 
(-60:1) circle(2pt) 
(60:1)++(0:2) circle(2pt)
(0:4) circle(2pt)
;

\end{tikzpicture}
\end{center}
\caption{Cyclically ordered 6-tuples of $o$'s and $x$'s for type-1 domain.}\label{f5}
\end{figure}
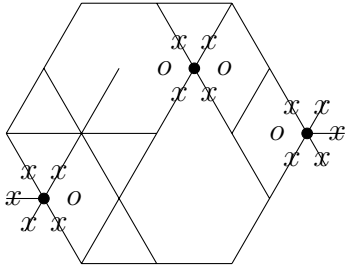

If for some edge $e$ of $U$ both endpoints of $e$ are not in $U$ (including the case of an endpoint in $\pa V_k$), then we can add $e$ to $L$ thus decreasing the number of edges in $U$ and leaving the number of vertices in $U$ unchanged, and so strictly increasing $\mu_1$. So we may assume that for every edge $e$ of $V_k$, if both endpoints of $e$ are not in $U$ then $e$ is not in $U$. If follows that for each triangle of $V_k$, if two of its sides are not in $U$ (including the case of a side contained in $\pa V_k$), then also its third side is not in $U$. This implies that all boundary angles of a connected component of $U$ which is not a single triangle, are at least $120^\circ$.

If for some vertex $v$ of $V_k$  the corresponding 6-tuple includes 1, 2, or 3 $x$'s
then we may remove the corresponding 1-simplices from $L$.
This adds 1, 2, or  3 edges  to $U$, and at least one vertex, namely $v$ itself. 
So if $a$ is the number of additional vertices and $b$ is the number of additional edges then $\frac{a}{b}\geq\frac{1}{3}$. From Lemmas \pr{01} and \pr{p}(1) it then follows that $\mu_1$ strictly increases. 
We may thus assume that for each vertex of $V_k$ the corresponding 6-tuple includes either 0 or at least 4 $x$'s. This is also true for the vertices of $\pa V_k$, since as we have noted above, each such vertex has at least 4 $x$'s. It follows that all boundary angles of the connected components of $U$ are at most $180^\circ$.

We thus need to establish our claim when
all boundary angles of the connected components of $U$ which are not single triangles, are either $120^\circ$ or $180^\circ$. It follows, as before, that each connected component of $U$ is either a single triangle or a hexagon. The single triangles have no vertices and no edges and so they do not contribute to $\mu_1(U)$. Thus by Lemma \pr{p}(1) it is enough to prove our claim for each of the hexagons.

Look at one such hexagon, and denote it again by $U$. We construct a larger hexagon $U'$ in $V_k$, as in the proof of Proposition \pr{u}, by shifting $U$ in $V_k$ if needed, using Lemma \pr{hx2}, and then adding a full row of triangles along the longest edge of $U$, whose length we denote $t$. This process adds
 $t-1$ vertices and $3t-2$ edges, as seen in Figure \pr{f4} where $t=4$. 
By Lemma \pr{hx} we have  
$\mu_1=\frac{f_2-b+2}{3f_2-b}$. 
To use Lemma \pr{p}(1) for showing that $\mu_1(U')>\mu_1(U)$ we need $\frac{t-1}{3t-2} > \frac{f_2-b+2}{3f_2-b}$, i.e.\ $f_2<(b-3)(2t-1)+1$.
The argument proceeds from here as in the proof of 
Proposition \pr{u}, this time leading to 
$\frac{\sqrt{3}}{\pi} < 1-\frac{6}{b}$, which holds for $b\geq 14$.

\end{pf}

We proceed to prove Proposition \pr{m}(1),
namely, that given $\epsilon>0$, then for sufficiently large $k$, every subcomplex $T$ of $\s S$ has $\mu_1(T)>\frac{1}{3}-\epsilon$.
Since $\mu_1$ does not involve 2-simplices, we may replace any $T$ with its 1-skeleton, and so it is enough to 
prove our claim for 1-dimensional subcomplexes, i.e.\  graphs. Furthermore, it is enough to prove the claim for pure 1-dimensional subcomplexes, i.e.\  graphs
with no isolated vertices, since given any graph $T$ we may remove its isolated vertices, which
may only reduce $\mu_1(T)$. Now,
given $\epsilon > 0$, by Proposition \pr{vk} we have $\mu_1(V_k)>\frac{1}{3}-\epsilon$ for sufficiently large $k$, and we show that every subgraph $T$ of $\s S$ with no isolated vertices satisfies $\mu_1(T)>\frac{1}{3}-\epsilon$. 
 For fixed 2-simplex $\D$ of $S$ we look at the corresponding $V_k$. If $T \cap V_k \neq \varnothing$ we claim that 
$\mu_1(T \cap V_k) > \frac{1}{3}-\epsilon$. 
Indeed this holds if $T \cap V_k$ includes all edges of $V_k$, since then $\mu_1(T \cap V_k)=\mu_1(V_k)$. Otherwise look at $U=V_k-(T\cap V_k)$, which is a type-1 domain in $V_k$ that includes edges and so $\mu_1(U)$ is defined. 
By Proposition \pr{u1} we have $\mu_1(U) < \mu_1(V_k)$, and so by Lemma \pr{p}(2) we have 
$\mu_1(T \cap V_k) = \mu_1(V_k-U) > \mu_1(V_k) > \frac{1}{3}-\epsilon$. If we look at the union of  $T \cap V_k$ over
all 2-simplices of $S$ then by Lemma \pr{p}(1) we still have $\mu_1>\frac{1}{3}-\epsilon$. 

This completes the computation for the intersection of $T$ with the interiors of the 
2-simplices of $S$, but we are still missing $T \cap \s S^1$ in case it is nonempty, (recall $S^1$ is the 1-skeleton of $S$). We claim that for $k\geq 1$, every subgraph of $\s S^1$ has 
$\mu_1>\frac{1}{2}>\frac{1}{3}-\epsilon$, which would complete our proof using Lemma \pr{p}(1).
Denoting the graph ${\mathfrak{s}}^{k-1}S^1$ by $P$, we have
$\s S^1={\mathfrak{s}}^1 P$. 
Note that for every edge $e$ of ${\mathfrak{s}}^1 P$, one endpoint of $e$ is an original vertex of $P$, and the other endpoint is a vertex subdividing an edge of $P$, and so its degree in ${\mathfrak{s}}^1 P$ is 2.
Now let 
$G\su {\mathfrak{s}}^1 P$ be a subgraph, and we may assume as before that $G$ has no isolated vertices. Let $v_1,\dots,v_m$ be those vertices of $G$ which are vertices of $P$, and let $d_i>0$ be the degree of $v_i$ in $G$. Let $H_i$ be $v_i$ together with its $d_i$ adjacent edges in $G$, and ``half'' of the vertex on the other side of each such edge. 
Then the $H_i$ are disjoint, and $H_i$ includes $1+\frac{d_i}{2}$ vertices and $d_i$ edges, and so 
$\mu_1(H_i)=\frac{1+\frac{d_i}{2}}{d_i} > \frac{1}{2}$. By Lemma \pr{p}(1) the combined contribution of all $H_i$ produces $\mu_1> \frac{1}{2}$. All that is now missing is perhaps some number of half vertices, which may only increase $\mu_1$,
so we finally get 
$\mu_1(G)> \frac{1}{2}$, completing the proof of Proposition \pr{m}(1).

Now that both parts of Proposition \pr{m} are proved, we have completed the proof of Theorem \pr{t1}(1).

\section{Topological non-embeddability}\label{ne}

In this concluding section we prove Theorem \pr{t1}(2) stating that for a probability triple $\p$ satisfying 
$\lim_{n\to\infty} np_0p_1^3p_2^2 = 0$, the torus is a.a.s.\ not topologically embeddable into $Y\in(\Omega_n,\pp)$. 
We will need the following lemma.

\begin{lemma}\label{l1}
Let $F$ be a closed surface and $Y$ a 2-dimensional simplicial complex. If $h:F \to Y$ is a topological embedding then $h(F)$ is a subcomplex of $Y$. It follows that there is a triangulation $S$ of $F$ for which
the embedding $h$ is simplicial.
\end{lemma}

\begin{pf}
We first show that if $V$ is the interior of a 2-simplex in $Y$, then either $h(F) \supseteq V$ or $h(F)\cap V=\varnothing$. This we do by showing that $h(F)\cap V$ is closed and open in $V$. Since $h(F)$ is compact, it is closed in $Y$, and so $h(F)\cap V$ is closed in $V$. On the other hand, since $Y$ is a 2-dimensional complex, $V$ is open in $Y$ and so $h^{-1}(V)$ is open in $F$. By invariance of domain (Corollary 19.9 in \ct{b}), $h|_{h^{-1}(V)}:h^{-1}(V)\to V$ is an open mapping, and so its image, $h(F)\cap V$, is open in $V$. 

Now let $G$ be the union of those 2-simplices in $Y$ whose interiors are contained in $h(F)$. Then $G$ is a subcomplex of $Y$ and we show $h(F)=G$. Since $h(F)$ is closed in $Y$ and is disjoint from the interiors of all other 2-simplices of $Y$, we have $G \su h(F) \su G\cup Y^1$, where $Y^1$ is the 1-skeleton of $Y$.  
Let $U=F-h^{-1}(G)$ then $U$ is open in $F$ and $h(U) \su Y^1$. Since $Y^1$ may be embedded into some surface as a subset with empty interior, by invariance of domain again we must have $U=\varnothing$ and so $h(F)=G$.
\end{pf}

By Lemma \pr{l1}, the torus is \emph{topologically} embeddable into a given 2-dimensional simplicial complex $Y$ iff there exists a triangulation $S$ of the torus that is \emph{simplicially} embeddable into $Y$.
If $S$ is a triangulation of the torus with $f_i$ simplices of dimension $i$, then there are $n(n-1)\cdots(n-f_0+1)\leq n^{f_0}$ injections of its vertex set into $[n]$. Given such an injection, the probability that it determines a simplicial embedding of $S$ into $Y\in(\Omega_n,\pp)$ is precisely $p_0^{f_0}p_1^{f_1}p_2^{f_2}$.
It follows that the probability that $S$ is simplicially embeddable into $Y$ is at most 
$n^{f_0} p_0^{f_0}p_1^{f_1}p_2^{f_2}$. 
Every triangulation of the torus satisfies $f_0-f_1+f_2 =0$ and $3f_2=2f_1$, from which it follows that 
$f_1=3f_0$ and $f_2=2f_0$. 
Furthermore, it follows from  \ct{g} that there is a constant $c$ such that the number of non-equivalent triangulations of the torus with $f_0$ vertices is bounded by $c^{f_0}$.
So, the probability that there is a triangulation $S$ of the torus with $f_0$ vertices which
is simplicially embeddable into $Y$ is at most $\big(cnp_0p_1^3p_2^2\big)^{f_0}=u_n^{f_0}$ where $u_n=cnp_0p_1^3p_2^2\longrightarrow 0$ by our assumption.
Finally, summing over all possible $f_0$, the probability that there is some triangulation of the torus
which is simplicially embeddable into $Y\in(\Omega_n,\pp)$
 is bounded by 
$$\sum_{f_0=1}^\infty u_n^{f_0}=
\frac{u_n}{1-u_n} \ \mathop{\longrightarrow}_{n \to \infty} \ 0.$$

%%%%%%%
\end{document}